\documentclass[12pt, letterpaper]{article}
\usepackage[T1]{fontenc}
\usepackage[latin9]{inputenc}
\usepackage{amsmath}
\usepackage{amssymb}\usepackage{amsfonts}
\usepackage[hidelinks]{hyperref}
\newcommand{\Z}{\mathbb{Z}}

\begin{document}
\begin{flushright}
\par\end{flushright}\bigskip{}

\begin{center}
\textbf{\huge{Integer Powers Expressed as Nested Sums of Lower Powers modulo 2, 3}}{\huge{}}
\\ \vspace{0.5cm} {\huge{ }}\large{\textsc{nikita gurin}}{\huge{}}\\

\par\end{center}{\huge \par}
\begin{abstract}In this paper, we explore identities that allow for representation of positive integers raised to positive integral powers as sums of nested sums of smaller positive integral powers. We begin by establishing the base identity involving consecutive descending powers, which we then employ to construct a simple proof of Fermat's Little Theorem. After that, the focus shifts to deriving and rigorously proving a more intricate identity that represents integer powers as nested sums of descending powers of the same parity, i.e., descending powers modulo 2. The discussion is concluded by stating a prominent identity that allows for the representation of integer powers by lower powers modulo 3. 
\end{abstract} 
\vspace{1cm}
\begin{center}{\textsc{1. Introduction}}\end{center}
Our starting point is the identity that describes how any positive integer $x$ raised to a positive integer $n$ can be expressed as the sum of consecutive powers  from 1 to $n-1$. Representing positive integers $x$ in such a manner naturally extends to all integers, once we observe that $(-x)^n=-|x^n|$ for odd $n$'s and $(-x)^n=|x^n|$ for even $n$'s, where $|\cdot|$ is the absolute value. The following will serve as the base for all the subsequent and more complicated results in this work.
\begin{equation}x^n=x+\sum_{k\,=\,1}^{x-1}{\Big((k+1)^n-(k^n+1)\Big).}\end{equation}
Let us verify that this identity is indeed true for all integer powers.
\\
\\
\textbf{Lemma 1.1} \textit{If a positive integer $x$ is raised to some $n$-th power, where $n \in \mathbb{Z^+},$ then $x^n$ can be expressed using $(1)$.}
\\
\\
\textit{Proof.} We will show that the identity holds via the Principle of Mathematical Induction. Let $x,n \in \mathbb{Z^+}$ and consider the expression
 \[x+\sum_{k\,=\,1}^{x-1}{\Big((k+1)^n-(k^n+1)\Big).}\]
 Let $S(i)$ be the statement that for each $i \leq x-1 \in \mathbb{Z^+}$,
 \[x+\sum_{k\,=\,i}^{x-1}{\Big((k+1)^n-(k^n+1)\Big)}+i^n-i=x+\sum_{k\,=\,1}^{x-1}{\Big((k+1)^n-(k^n+1)\Big)}.\]
Base case: For $i=1$, we have
 \[x+\sum_{k\,=\,i}^{x-1}{\Big((k+1)^n-(k^n+1)\Big)}+i^n-i=x+\sum_{k\,=\,1}^{x-1}{\Big((k+1)^n-(k^n+1)\Big)}+1^n-1,\]
so $S(1)$ is true.
\\ \\ Induction step: Assume $S(i)$ is true. Then, 
 \[x+\sum_{k\,=\,i}^{x-1}{\Big((k+1)^n-(k^n+1)\Big)}+i^n-i=x+\sum_{k\,=\,1}^{x-1}{\Big((k+1)^n-(k^n+1)\Big)}.\]
Consequently, we compute the term of the summation on the left hand side at $k=i$ to get 
 \[x+\sum_{k\,=\,i+1}^{x-1}{\Big((k+1)^n-(k^n+1)\Big)}+\Big((i+1)^n-(i^n+1)\Big)+i^n-i=x+\sum_{k\,=\,1}^{x-1}{\Big((k+1)^n-(k^n+1)\Big)},\]
which allows us to find that
\[x+\sum_{k\,=\,i+1}^{x-1}{\Big((k+1)^n-(k^n+1)\Big)}+(i+1)^n-(i+1)=x+\sum_{k\,=\,1}^{x-1}{\Big((k+1)^n-(k^n+1)\Big)},\]
and thus $S(i+1)$ is true. Therefore, $S(i)$ is true for all $i \leq x-1 \in \mathbb{Z^+}$. 
\\ \\ Now, by the above result, note that
\[x+\sum_{k\,=\,x-1}^{x-1}{\Big((k+1)^n-(k^n+1)\Big)}+(x-1)^n-(x-1)=x+\sum_{k\,=\,1}^{x-1}{\Big((k+1)^n-(k^n+1)\Big)}.\]
Direct evaluation of the summation on the left at $k=x-1$ yields
\[x+(x-1+1)^n-\big((x-1)^n+1\big)+(x-1)^n-(x-1)=x+\sum_{k\,=\,1}^{x-1}{\Big((k+1)^n-(k^n+1)\Big)},\]
and simplifying the expression on the right finally reveals 
\begin{center}\vspace{-0.5cm}\[x^n=x+\sum_{k\,=\,1}^{x-1}{\Big((k+1)^n-(k^n+1)\Big)}, \ \textrm{as desired.}\] \hfill$\square$\end{center}
One powerful consequence of Lemma 1.1 is an easy approach to justification of the well-known Fermat's Little Theorem \hyperlink{ref1}{$[1]$}. In what follows, we utilize (1) to prove the theorem.    
\\ \\ \\ \\ \textbf{Theorem 1.1} (Fermat): \textit{$x^p \equiv x$ \emph{mod} $p$ \ for all $\ x,p \in \mathbb{Z^+}$ with \hspace{0.005cm} $p \in$ prime.}
\\ \\ \textit{Proof}. Let $x,p \in \mathbb{Z^+}$ such that $p$ is prime, and, by Lemma 1.1, assume
\[x^p=x+\sum_{k\,=\,1}^{x-1}{\Big((k+1)^p-k^p-1\Big).}\]
Observe that \[(k+1)^p-k^p-1=\sum_{i\,=\,1}^{p-1}{\frac{p!}{i!(p-i)!}k^i} \ \ \textrm{for each $k \leq x-1$}.\] Since $p$  is prime, gcd$(p,i!)=1$ and gcd$(p,(p-i)!)=1$ for all $i \in \mathbb{Z^+}$ such that $i \leq p-1$. Furthermore, since \[\frac{p!}{i!(p-i)!} = \frac{p \cdot (p-1)!}{i!(p-i)!} \in \mathbb{Z^+}\] as well as gcd$(p,i!)=$ gcd$(p,(p-i)!)=1$, by Euclid's Lemma \hyperlink{ref2}{$[2]$}, \[i!(p-i)! \mid (p-1)! \ \textrm{for every $i = 1,2,...,p-1$}.\] Let $\frac{(p-1)!}{i!(p-i)!}=p_i$ with $p_i \in \mathbb{Z^+}$ and $1 \leq i \leq p-1$. Then, we have 
\[x^p-x=\sum_{k\,=\,1}^{x-1}{\Big((k+1)^p-k^p-1\Big)}=p\cdot \sum_{k\,=\,1}^{x-1}{\sum_{i\,=\,1}^{p-1}{p_ik^i}},\] 
from which we deduce that $p \mid x^p-x$, and thus $x^p \equiv x$ mod $p$. \hfill $\square$ 
\\
\\
\indent After some experimenting, it turns out that a version of (1) is useful for representing integers $x^n$, namely
\begin{equation}x^n=\sum_{k\,=\,1}^{x}{k^n-(k-1)^n}.\end{equation}
\textbf{Lemma 1.2} \textit{If a positive integer $x$ is raised to some $n$-th power with $n \in \mathbb{Z^+}$, then $x^n$ can be expressed using $(2)$.}
\\
\\
\textit{Proof.} By Lemma 1.1, let $x,n \in \mathbb{Z^+}$ such that $x^n=x+\sum_{k\,=\,1}^{x-1}{\Big((k+1)^n-k^n-1\Big)}$.
Then, shifting the index of the summation by $1$ allows us to obtain
\[x^n=x+\sum_{k\,=\,2}^{x}{\Big(k^n-(k-1)^n-1\Big).}\]
Now, observe that at $k=1$, $k^n-(k-1)^n-1=1^n-(1-1)^n-1=0$, which means that $\sum_{k\,=\,2}{k^n-(k-1)^n-1}=\sum_{k\,=\,1}{k^n-(k-1)^n-1}$. Furthermore, $\sum_{k\,=\,1}^{x}{(-1)}=-x$, so the above results in
\begin{center}\[x^n=\sum_{k\,=\,1}^{x}{\Big(k^n-(k-1)^n\Big)}.\] \hfill $\square$\end{center}
\begin{center}{\textsc{2. Even and Odd Representations}}\end{center}
\hspace{0.45cm} Now, we transition to the consideration of a fairly intriguing identity that allows for the representation of any even power as sums of descending even powers and, likewise, any odd power as sums of descending odd powers. We state that identity as a theorem below and devote the rest of this section to proving it.
\\
\\
\textbf{Theorem 2.1} \textit{For all positive integers $x$ and $n$,}
\begin{equation}x^n=x+\sum_{s\,=\,1}^{m}{\frac{2n!}{(2s)!(n-2s)!}\sum_{k\,=\,1}^{x-1}{\sum_{j\,=\,1}^{k}{j^{n-2s}}}},\end{equation}
\textit{where} $m= \begin{cases} 
      \frac{n}{2}&\textrm{if} \ n \textrm{ is even} \\
      \frac{n-1}{2}&\textrm{if} \ n \textrm{ is odd} 
\end{cases}$.
\\
\\
\\ Before moving forward, let's take a look at \emph{how} the identity (3) works via the example that follows.
\\
\\
\textbf{Example} \textit{Let $x=5$ and $n=7$. Then, $5^7=78125$. At the same time, consider} \[5^7 = 5+\sum_{s\,=\,1}^{\frac{7-1}{2}}{\frac{2\cdot 7!}{(2s)!(7-2s)!}\sum_{k\,=\,1}^{5-1}{\sum_{j\,=\,1}^{k}{j^{7-2s}}}},\]
\textit{which can be restated as}
\[5^7 = 5+\frac{2\cdot 7!}{2!(7-2)!}\sum_{k\,=\,1}^{5-1}{\sum_{j\,=\,1}^{k}{j^{7-2}}}+\frac{2\cdot 7!}{4!(7-4)!}\sum_{k\,=\,1}^{5-1}{\sum_{j\,=\,1}^{k}{j^{7-4}}}+\frac{2\cdot 7!}{6!(7-6)!}\sum_{k\,=\,1}^{5-1}{\sum_{j\,=\,1}^{k}{j^{7-6}}},\]
\textit{so that}
\[5^7 = 5+2\cdot \frac{7!}{2!5!}\sum_{k\,=\,1}^{4}{\sum_{j\,=\,1}^{k}{j^5}}+2\cdot \frac{7!}{4!3!}\sum_{k\,=\,1}^{4}{\sum_{j\,=\,1}^{k}{j^3}}+2\cdot \frac{7!}{6!1!}\sum_{k\,=\,1}^{4}{\sum_{j\,=\,1}^{k}{j^1}},\]
\textit{and finally}
\[5^7=5+67620+10220+280=78125.\]
Now, we present a justification of \emph{why} the identity (3) works. 
\\
\\
\textit{Proof.} Let $x,n \in \mathbb{Z^+}$ and, by Lemma 1.1, assume
\[x^n=x+\sum_{k\,=\,1}^{x-1}{\Big((k+1)^n-(k^n+1)\Big)}.\]
Then, expanding inside the summation results in
\[x^n=x+\sum_{k\,=\,1}^{x-1}{\Bigg(\frac{n!}{1!(n-1)!}k^{n-1}+\frac{n!}{2!(n-2)!}k^{n-2}+ \ \dots \ +\frac{n!}{(n-1)!1!}k\Bigg)}.\]
Next, applying (2) to each $k^{n-\varepsilon}$ term with $1 \leq \varepsilon \leq n-1$ produces
\begin{align}
x^n= x + \sum_{k\,=\,1}^{x-1}\Bigg(\ & \frac{n!}{1!(n-1)!}\bigg(\sum_{j\,=\,1}^{k}{j^{n-1}-(j-1)^{n-1}\bigg)+\frac{n!}{2!(n-2)!}\bigg(\sum_{j\,=\,1}^{k}{j^{n-2}-(j-1)^{n-2}}\bigg)} + \notag \\
&+\ \dots \ + \frac{n!}{(n-2)!2!}\bigg(\sum_{j\,=\,1}^{k}{j^2-(j-1)^2}\bigg) + \frac{n!}{(n-1)!1!}\sum_{j\,=\,1}^{k}{j^0}\Bigg).
\notag \end{align}
Observe that there are exactly $n-1$ summations $\sum_{j\,=\,1}^{k}$ inside $\sum_{k\,=\,1}^{x-1}$. The following is obtained by expanding inside each $\sum_{j\,=\,1}^{k}$.
\begin{align}
x^n=x+\sum_{k\,=\,1}^{x-1}\Bigg(\ & \frac{n!}{1!(n-1)!}\sum_{j\,=\,1}^{k}{\Bigg(\frac{(n-1)!}{1!(n-2)!}j^{n-2}-\frac{(n-1)!}{2!(n-3)!}j^{n-3}+\ \dots \ +\frac{(n-1)!}{(n-2)!1!}j^0\Bigg)} \notag\\
&+ \frac{n!}{2!(n-2)!}\sum_{j\,=\,1}^{k}{\Bigg(\frac{(n-2)!}{1!(n-3)!}j^{n-3}-\frac{(n-2)!}{2!(n-4)!}j^{n-4}+\ \dots \ -\frac{(n-2)!}{(n-3)!1!}j^0\Bigg)} \notag \\
&+ \ \dots \ + \notag \\
&+ \frac{n!}{\varepsilon!(n-\varepsilon)!} \sum_{j\,=\,1}^{k}{\Bigg(\frac{(n-\varepsilon)!}{1!(n-\varepsilon-1)!}j^{n-\varepsilon-1} + \ \dots \ +(-1)^{\lambda+1}\frac{(n-\varepsilon)!}{\lambda!(n-\varepsilon-\lambda)!}j^{n-\varepsilon-\lambda}} \notag \\ 
& + \ \dots \ + \ (-1)^{\varepsilon+1}\frac{(n-\varepsilon)!}{(n-\varepsilon-1)!1!}j^0\Bigg) \notag \\
& + \ \dots \ + \notag \\
& + \frac{n!}{(n-1)!1!}\sum_{j\,=\,1}^{k}{\frac{1!}{1!0!}j^0}\Bigg)
\end{align}
with each $\sum_{j\,=\,1}^{k}$ having exactly $n-\varepsilon$ summands of the form $(-1)^{\lambda+1}\frac{(n-\varepsilon)!}{\lambda!(n-\varepsilon-\lambda)!}j^{n-\varepsilon-\lambda}$ in the definition of its terms, where $\varepsilon$ is an increasing index of a given $\sum_{j\,=\,1}^{k}$ such that $\varepsilon=1,2, \dots, n-1$. Note that each summation has one less summand $(-1)^{\lambda+1}\frac{(n-\varepsilon)!}{\lambda!(n-\varepsilon-\lambda)!}j^{n-\varepsilon-\lambda}$ than the summation that immediately precedes it in the above top-to-bottom arrangement of (4) because $\varepsilon$ increments by 1 with each new summation as we go from top to bottom.
\\ \\ Let $r=\varepsilon+\lambda$ so that $2 \leq r \leq n$, where $\lambda$ is a decreasing index of the position of a given summand $(-1)^{\lambda+1}\frac{(n-\varepsilon)!}{\lambda!(n-r)!}j^{n-r}$ in each $\sum_{j\,=\,1}^{k}$ with $\lambda=r-\varepsilon,\dots,2,1$. Observe that $\lambda$ decrements by 1 with each new summation as we go from top to bottom. Hence, $r$ stays unchanged from one summation to the next in (4), since all increases in $\varepsilon$ are canceled by equal decreases in $\lambda$. For instance, in the first sum (the uppermost in our arrangement), $\varepsilon=1$, so $\lambda=r-1$; then, in the next sum below, $\varepsilon=2$, so $\lambda=r-2$, such that they still add up to $r$, etc.
\\ \\ Note that for a given $r$, each summand of the form $(-1)^{\lambda+1}\frac{(n-\varepsilon)!}{\lambda!(n-r)!}j^{n-r}$, we call it a $j^{n-r}$-like term, appears either once or zero times in each $\sum_{j\,=\,1}^{k}$. Since $\lambda \geq 1$, the lowest value in our top-to-bottom arrangement (hence the last value) that $r$ is represented by is $\varepsilon_{\lambda\,=\,1}+1$ for some $\varepsilon_{\lambda\,=\,1}=1,2, \dots, n-1$, so the only $\sum_{j\,=\,1}^{k}$'s that have $j^{n-r}$-like terms $(-1)^{\lambda+1}\frac{(n-\varepsilon)!}{\lambda!(n-r)!}j^{n-r}$ for a given $r$ are those that are ordered from $1$ to $\varepsilon_{\lambda\,=\,1}$ from top to bottom, and there are $\varepsilon_{\lambda\,=\,1}=r-1$ of them. Thus, for each $r$, exactly $r-1$ summations $\sum_{j\,=\,1}^{k}$ contain exactly one $j^{n-r}$-like term. 
\\
\\
For example, given $r=3$, $j^{n-3}$-like terms appear in exactly $3-1=2$ summations of (4), namely
\[\frac{n!}{1!(n-1)!}\sum_{j\,=\,1}^{k}{\Bigg(\frac{(n-1)!}{1!(n-2)!}j^{n-2}\underbrace{-\frac{(n-1)!}{2!(n-3)!}j^{n-3}}_{\textrm{here}}+\ \dots \ +\frac{(n-1)!}{(n-2)!1!}j^0\Bigg)} \ \textrm{with $\varepsilon=1$, $\lambda=2$}\]
and
\[\frac{n!}{2!(n-2)!}\sum_{j\,=\,1}^{k}{\Bigg(\underbrace{\frac{(n-2)!}{1!(n-3)!}j^{n-3}}_{\textrm{here}}-\frac{(n-2)!}{2!(n-4)!}j^{n-4}+\ \dots \ -\frac{(n-2)!}{(n-3)!1!}j^0\Bigg)} \ \textrm{with $\varepsilon_{\lambda\,=\,1}=2$, $\lambda=1$}.\]
Now, the idea is to take all $j^{n-r}$-like terms and reorganize them in a more useful order. To do this, we take all  summands $(-1)^{r-\varepsilon+1}\frac{(n-\varepsilon)!}{(r-\varepsilon)!(n-r)!}j^{n-r}$ for each $r=2,3,\dots,n$ (remember there are always $r-1$ of such summands), going from first down to each $\varepsilon_{\lambda\,=\,1}$-th summation, and group them together, thereby separating them from the rest of the summations. Below is an example of grouping $j^{n-4}$-like terms together. 
\begin{align}
x^n=x+\sum_{k\,=\,1}^{x-1}\Bigg(\ & \frac{n!}{1!(n-1)!}\sum_{j\,=\,1}^{k}{\Bigg(\dots\Bigg)} + \frac{n!}{1!(n-1)!}\frac{(n-1)!}{3!(n-4)!}j^{n-4} \notag\\
&+ \frac{n!}{2!(n-2)!}\sum_{j\,=\,1}^{k}{\Bigg(\dots\Bigg)} - \frac{n!}{2!(n-2)!}\frac{(n-2)!}{2!(n-4)!}j^{n-4} \notag \\
&+ \frac{n!}{3!(n-3)!}\sum_{j\,=\,1}^{k}{\Bigg(\dots\Bigg)} + \frac{n!}{3!(n-3)!}\frac{(n-3)!}{1!(n-4)!}j^{n-4} \notag \\
& \hspace{-0.113cm}\left.
\begin{aligned}
&+ \ \dots \ + \\
&+ \frac{n!}{(n-1)!1!}\sum_{j\,=\,1}^{k}{\frac{1!}{1!0!}j^0}\Bigg)                  
\end{aligned}
\right\} \ \text{the rest of the summations without $j^{n-4}$-like terms} \notag
\end{align}
becomes
\begin{align}
x^n=x+\sum_{k\,=\,1}^{x-1}\Bigg(\ & \frac{n!}{1!(n-1)!}\frac{(n-1)!}{3!(n-4)!}j^{n-4} - \frac{n!}{2!(n-2)!}\frac{(n-2)!}{2!(n-4)!}j^{n-4} + \frac{n!}{3!(n-3)!}\frac{(n-3)!}{1!(n-4)!}j^{n-4} \notag \\
& +\underbrace{\sum_{j\,=\,1}^{k}{\Bigg(\dots\Bigg)} }_\text{the rest of the summations without $j^{n-4}$-like terms} \Bigg). \notag
\end{align}
As we perform this grouping action to $j^{n-r}$-like terms for each $r$, we produce
\begin{align}
x^n=x+\sum_{r\,=\,2}^{n}\sum_{k\,=\,1}^{x-1}\sum_{j\,=\,1}^{k}j^{n-r}\Bigg( \ & (-1)^{r}\frac{n!}{1!(n-1)!}\frac{(n-1)!}{(r-1)!(n-r)!}+(-1)^{r-1}\frac{n!}{2!(n-2)!}\frac{(n-2)!}{(r-2)!(n-r)!} \notag \\
&+ \ \dots \ + (-1)^{r-i+1}\frac{n!}{i!(n-i)!}\frac{(n-i)!}{(r-i)!(n-r)!} \notag \\
&+ \ \dots \ + (-1)^{2}\frac{n!}{\varepsilon_{\lambda\,=\,1}!(n-\varepsilon_{\lambda\,=\,1})!}\frac{(n-r+1)!}{1!(n-r)!}\Bigg) \notag
\end{align}
(note that $\varepsilon_{\lambda\,=\,1} = r-1$), in which we treat the resulting summation of fractions as a coefficient of the nested summation $\sum_{k\,=\,1}^{x-1}{\sum_{j\,=\,1}^{k}{j^{n-r}}}$ and factor it out as below.
\[x^n=x+\sum_{r\,=\,2}^{n}\frac{n!}{(n-r)!}\Bigg(\sum_{i\,=\,1}^{r-1}{\frac{(-1)^{r-i+1}}{i!(r-i)!}\Bigg)\cdot\sum_{k\,=\,1}^{x-1}{\sum_{j\,=\,1}^{k}{j^{n-r}}}}.\]
Consider the nested summation on the right hand side of the above for a particular $r$
\begin{equation}\frac{n!}{(n-r)!}\sum_{i\,=\,1}^{r-1}{\frac{(-1)^{r-i+1}}{i!(r-i)!}\cdot\sum_{k\,=\,1}^{x-1}{\sum_{j\,=\,1}^{k}{j^{n-r}}}}.\end{equation}
Let us examine (5) further by describing what this expression evaluates to for different $r$'s. We proceed by cases.
\\ \\ 
Case 1: Suppose $r$ is odd. Then, $r \geq 3$ and $r=2s+1$ for some $s \in \mathbb{Z^+}$. Consequently, (5) turns into 
\[\frac{n!}{(n-2s-1)!}\sum_{i\,=\,1}^{2s}{\frac{(-1)^{2s+2-i}}{i!(2s+1-i)!}\sum_{k\,=\,1}^{x-1}{\sum_{j\,=\,1}^{k}{j^{n-2s-1}}}}.\]
Using the Principle of Mathematical Induction, we will now show that $\sum_{i\,=\,1}^{2s}{\frac{(-1)^{2s+2-i}}{i!(2s+1-i)!}}=0$ for $s \in \mathbb{Z^+}$, and hence, (5) always evaluates to $0$ when $r$ is odd.
\[\textrm{Let} \ P(t): \sum_{i\,=\,1}^{2s}{\frac{(-1)^{2s+2-i}}{i!(2s+1-i)!}}=\sum_{i\,=\,1+t}^{2s-t}{\frac{(-1)^{2s+2-i}}{i!(2s+1-i)!}\ \textrm{for each $t \in \mathbb{Z^+}$}}.\] The motivation behind the lines of reasoning that follow is that for any given $s$ and for sufficiently a large $t$, it is true that $2s-t < 1+t$, so $\sum_{i\,=\,1+t}^{2s-t}{\frac{(-1)^{2s+2-i}}{i!(2s+1-i)!}}$ is a sum whose upper bound is smaller than its lower bound, which, by definition, is an empty sum and is equal to $0$. Hence, $\sum_{i\,=\,1}^{2s}{\frac{(-1)^{2s+2-i}}{i!(2s+1-i)!}}$ must also be equal to 0.
\\
\\
Base step: Observe that, by taking the first and last terms, i.e., terms at $i=1$ and $i=2s$, respectively, out of the summation $\sum_{i\,=\,1}^{2s}{\frac{(-1)^{2s+2-i}}{i!(2s+1-i)!}}$, we get
\begin{align}
\sum_{i\,=\,1}^{2s}{\frac{(-1)^{2s+2-i}}{i!(2s+1-i)!}}&=\frac{(-1)^{2s+2-1}}{1!(2s+1-1)!}+\frac{(-1)^{2s+2-2s}}{(2s)!(2s+1-2s)!}+\sum_{i\,=\,2}^{2s-1}{\frac{(-1)^{2s+2-i}}{i!(2s+1-i)!}} \notag \\
&= \frac{-1}{1!(2s)!}+\frac{1}{(2s)!1!}+\sum_{i\,=\,2}^{2s-1}{\frac{(-1)^{2s+2-i}}{i!(2s+1-i)!}} \notag \\
&=\sum_{i\,=\,1+1}^{2s-1}{\frac{(-1)^{2s+2-i}}{i!(2s+1-i)!}}, \notag
\end{align}
so $P(1)$ is true.
\\
\\
Induction step: Assume $P(t)$ is true. Then,
\[\sum_{i\,=\,1}^{2s}{\frac{(-1)^{2s+2-i}}{i!(2s+1-i)!}}=\sum_{i\,=\,1+t}^{2s-t}{\frac{(-1)^{2s+2-i}}{i!(2s+1-i)!}}.\]
Note the fact that taking the first and last terms, i.e., terms at $i=1+t$ and $i=2s-t$, respectively, out of the summation $\sum_{i\,=\,1+t}^{2s-t}{\frac{(-1)^{2s+2-i}}{i!(2s+1-i)!}}$ yields
\begin{align}
\sum_{i\,=\,1+t}^{2s-t}{\frac{(-1)^{2s+2-i}}{i!(2s+1-i)!}}&=\frac{(-1)^{2s+2-(1+t)}}{(1+t)!(2s+1-(1+t))!}+\frac{(-1)^{2s+2-(2s-t)}}{(2s-t)!(2s+1-(2s-t))!}+\sum_{i\,=\,1+t+1}^{2s-t-1}{\frac{(-1)^{2s+2-i}}{i!(2s+1-i)!}} \notag \\
& =\frac{(-1)^{2s+1-t}}{(1+t)!(2s-t)!}+\frac{(-1)^{2+t}}{(2s-t)!(1+t)!}+\sum_{i\,=\,1+(t+1)}^{2s-(t+1)}{\frac{(-1)^{2s+2-i}}{i!(2s+1-i)!}}. \notag 
\end{align}
Here, we consider cases.
\\
\\
Case 1.a: Suppose $t$ is odd. Then, $(2s+1)-t=$ odd $-$ odd $=$ even, $s \in \Z^+$. Also, $2+t=$ even $+$ odd $=$ odd. This means that
\[\frac{(-1)^{2s+1-t}}{(1+t)!(2s-t)!}+\frac{(-1)^{2+t}}{(2s-t)!(1+t)!}=\frac{1}{(1+t)!(2s-t)!}+\frac{-1}{(2s-t)!(1+t)!}=0.\]
\\
Case 1.b: Suppose $t$ is even. Then, $(2s+1)-t=$ odd $-$ even $=$ odd. On the other hand, $2+t=$ even $+$ even $=$ even. This results in
\[\frac{(-1)^{2s+1-t}}{(1+t)!(2s-t)!}+\frac{(-1)^{2+t}}{(2s-t)!(1+t)!}=\frac{-1}{(1+t)!(2s-t)!}+\frac{1}{(2s-t)!(1+t)!}=0.\]
\\
In all cases, $\frac{(-1)^{2s+1-t}}{(1+t)!(2s-t)!}+\frac{(-1)^{2+t}}{(2s-t)!(1+t)!}=0$, and thus,
\[\sum_{i\,=\,1}^{2s}{\frac{(-1)^{2s+2-i}}{i!(2s+1-i)!}}\shortstack[pos]{\scriptsize{induction hypothesis}\\ \ \ $=$ \ }\sum_{i\,=\,1+t}^{2s-t}{\frac{(-1)^{2s+2-i}}{i!(2s+1-i)!}} \ =\sum_{i\,=\,1+(t+1)}^{2s-(t+1)}{\frac{(-1)^{2s+2-i}}{i!(2s+1-i)!}},\]
so $P(t+1)$ is true. Therefore, $P(t)$ is true for all $t \in \mathbb{Z^+}$. This implies that there always exists $t$ such that, for any given $s$, $2s-t < 1+t$. Taking advantage of this fact and picking a sufficiently large $t$ for any $s$ allows us to deduce that $\sum_{i\,=\,1+t}^{2s-t}{\frac{(-1)^{2s+2-i}}{i!(2s+1-i)!}}$ is an empty sum, so that
\[\sum_{i\,=\,1}^{2s}{\frac{(-1)^{2s+2-i}}{i!(2s+1-i)!}}=\sum_{i\,=\,1+t}^{2s-t}{\frac{(-1)^{2s+2-i}}{i!(2s+1-i)!}}=0,\]
which, in turn, yields 
\[\frac{n!}{(n-2s-1)!}\sum_{i\,=\,1}^{2s}{\frac{(-1)^{2s+2-i}}{i!(2s+1-i)!}\sum_{k\,=\,1}^{x-1}{\sum_{j\,=\,1}^{k}{j^{n-2s-1}}}}=\frac{n!}{(n-2s-1)!}\cdot0\cdot\sum_{k\,=\,1}^{x-1}{\sum_{j\,=\,1}^{k}{j^{n-2s-1}}}=0.\] Therefore, (5) is equal to 0 for all odd $r$'s.
\\
\\
\\
Case 2: Suppose $r$ is even. Then, $r \geq 2$ and $r=2s$ for some $s \in \mathbb{Z^+}$. Next, we rewrite (5) as 
\[\frac{n!}{(n-2s)!}\sum_{i\,=\,1}^{2s-1}{\frac{(-1)^{2s+1-i}}{i!(2s-i)!}\sum_{k\,=\,1}^{x-1}{\sum_{j\,=\,1}^{k}{j^{n-2s}}}}.\]
We proceed by induction to show that $\sum_{i\,=\,1}^{2s-1}{\frac{(-1)^{2s+1-i}}{i!(2s-i)!}}=\frac{2}{(2s)!}$ for $s \in \mathbb{Z^+}.$
\[\textrm{Let} \ Q(s): \sum_{i\,=\,1}^{2s-1}{\frac{(-1)^{2s+1-i}}{i!(2s-i)!}}=\frac{2}{(2s)!} \ \textrm{for each $s \in \mathbb{Z^+}$}.\]
Base step: at $s=1$, $\sum_{i\,=\,1}^{2s-1}{\frac{(-1)^{2s+1-i}}{i!(2s-i)!}}=\sum_{i\,=\,1}^{1}{\frac{(-1)^{3-i}}{i!(2-i)!}}=1=\frac{2}{(2\cdot1)!},$ so $Q(1)$ is true.
\\
\\
Induction step: Assume $Q(s)$ is true. Then,
\[\sum_{i\,=\,1}^{2s-1}{\frac{(-1)^{2s+1-i}}{i!(2s-i)!}}=\frac{2}{(2s)!}.\]
By the induction hypothesis and using the identity $\frac{2}{(2s)!}+\frac{2(1-(2s+1)(2s+2))}{(2s+2)!}=\frac{2}{(2s+2)!}$, we deduce 
\begin{equation}\sum_{i\,=\,1}^{2s-1}{\Bigg(\frac{(-1)^{2s+1-i}}{i!(2s-i)!}\Bigg)}+\frac{2(1-(2s+1)(2s+2))}{(2s+2)!}=\frac{2}{(2s+2)!}.\end{equation}
Now, consider the following identity.
\begin{equation}\frac{-1}{(2s)!2!}+\frac{1}{(2s+1)!1!}+\frac{2((2s+1)(2s+2)-1)}{(2s+2)!}=\sum_{i\,=\,1}^{2s-1}{\frac{(-1)^{1+i}\big((2s-i)^2+3(2s-i)+1\big)}{i!(2s+2-i)!}}.\end{equation}
Also, observe that
\begin{align}\sum_{i\,=\,1}^{2s-1}{\frac{(-1)^{1+i}\big((2s-i)^2+3(2s-i)+1\big)}{i!(2s+2-i)!}}&=\sum_{i\,=\,1}^{2s-1}{\frac{(-1)^{1+i}\big((2s+1-i)(2s+2-i)-1\big)}{i!(2s+2-i)!}} \notag
\\ &=\sum_{i\,=\,1}^{2s-1}{\Bigg(\frac{(-1)^{1+i}(2s+1-i)(2s+2-i)}{i!(2s+2-i)!}-\frac{(-1)^{1+i}}{i!(2s+2-i)!}\Bigg)}, \notag \end{align}
and so, by virtue of the facts that $(-1)^{i}=(-1)^{-i}$ for all $i \in \mathbb{Z^+}$ and $(-1)^1=(-1)^{2s+1}=$ \\ $=(-1)^{2s+3}$ for $s \in \mathbb{Z^+}$, we have
\begin{equation}\sum_{i\,=\,1}^{2s-1}{\frac{(-1)^{1+i}\big((2s-i)^2+3(2s-i)+1\big)}{i!(2s+2-i)!}}=\sum_{i\,=\,1}^{2s-1}{\Bigg(\frac{(-1)^{2s+1-i}}{i!(2s-i)!}\Bigg)}-\sum_{i\,=\,1}^{2s-1}{\Bigg(\frac{(-1)^{2s+3-i}}{i!(2s+2-i)!}\Bigg)}.\end{equation}
Combining the results of (7) and (8) produces
\[\frac{-1}{(2s)!2!}+\frac{1}{(2s+1)!1!}+\frac{2((2s+1)(2s+2)-1)}{(2s+2)!}=\sum_{i\,=\,1}^{2s-1}{\Bigg(\frac{(-1)^{2s+1-i}}{i!(2s-i)!}\Bigg)}-\sum_{i\,=\,1}^{2s-1}{\Bigg(\frac{(-1)^{2s+3-i}}{i!(2s+2-i)!}\Bigg)}.\]
Then, rearranging allows us to obtain
\[\sum_{i\,=\,1}^{2s-1}{\frac{(-1)^{2s+3-i}}{i!(2s+2-i)!}}+\frac{-1}{(2s)!2!}+\frac{1}{(2s+1)!1!}=\sum_{i\,=\,1}^{2s-1}{\frac{(-1)^{2s+1-i}}{i!(2s-i)!}}+\frac{2(1-(2s+1)(2s+2))}{(2s+2)!}.\]
Next, applying the result of (6) on the right-hand side of the above yields
\[\sum_{i\,=\,1}^{2s-1}{\frac{(-1)^{2s+3-i}}{i!(2s+2-i)!}}+\frac{-1}{(2s)!2!}+\frac{1}{(2s+1)!1!}=\frac{2}{(2s+2)!}.\]
Consequently, notice that $\frac{-1}{(2s)!2!}$ and $\frac{1}{(2s+1)!1!}$ precisely match the form of terms of $\sum_{i\,=\,1}^{2s-1}{\frac{(-1)^{2s+3-i}}{i!(2s+2-i)!}}$:
\begin{align}&\textrm{at $i=2s$}, \ \frac{(-1)^{2s+3-i}}{i!(2s+2-i)!}=\frac{(-1)^{2s+3-(2s)}}{(2s)!(2s+2-(2s))!}=\frac{-1}{(2s)!2!} \notag; 
\\ &\textrm{at $i=2s+1$}, \ \frac{(-1)^{2s+3-i}}{i!(2s+2-i)!}=\frac{(-1)^{2s+3-(2s+1)}}{(2s+1)!(2s+2-(2s+1))!}=\frac{1}{(2s+1)!1!} \notag; 
\end{align}
so we merge these two terms into the summation, which results in
\[\sum_{i\,=\,1}^{2s+1}{\frac{(-1)^{2s+3-i}}{i!(2s+2-i)!}}=\frac{2}{(2s+2)!},\]
or, expressed more explicitly,
\[\sum_{i\,=\,1}^{2(s+1)-1}{\frac{(-1)^{2(s+1)+1-i}}{i!(2(s+1)-i)!}}=\frac{2}{(2(s+1))!},\]
which shows that $Q(s+1)$ is true. Therefore,
\[\sum_{i\,=\,1}^{2s-1}{\frac{(-1)^{2s+1-i}}{i!(2s-i)!}}=\frac{2}{(2s)!} \ \textrm{for all $s \in \mathbb{Z^+}$},\]
which means that (5) becomes
\[\frac{2n!}{(2s)!(n-2s)!}\sum_{k\,=\,1}^{x-1}{\sum_{j\,=\,1}^{k}{j^{n-2s}}}\]
for all $r=2s$, $s \in \mathbb{Z^+}$, i.e., for all even $r$'s. Now, recall that 
\[x^n=x+\sum_{r\,=\,2}^{n}{\frac{n!}{(n-r)!}\sum_{i\,=\,1}^{r-1}{\frac{(-1)^{r-i+1}}{i!(r-i)!}\sum_{k\,=\,1}^{x-1}{\sum_{j\,=\,1}^{k}{j^{n-r}}}}},\]
which is the same as
\[x^n=x+\Bigg(\sum_{\textrm{odd $r$'s, $3 \leq r \leq n$}}+\sum_{\textrm{even $r$'s, $2 \leq r \leq n$}}\Bigg){\frac{n!}{(n-r)!}}\sum_{i\,=\,1}^{r-1}{\frac{(-1)^{r-i+1}}{i!(r-i)!}}\sum_{k\,=\,1}^{x-1}{\sum_{j\,=\,1}^{k}{j^{n-r}}}.\]
Since, by consideration of cases, $\sum_{i\,=\,1}^{r-1}{\frac{(-1)^{r-i+1}}{i!(r-i)!}} = \begin{cases} 
      0 & r=2s+1, \ s \in \mathbb{Z^+} \\
      \frac{2}{(2s)!} & r=2s, \ s \in \mathbb{Z^+} 
\end{cases},$ we have
\[x^n=x+\Bigg(0+\sum_{s\,=\,\frac{r}{2}, \ 1 \leq s \leq \frac{n}{2}}\frac{2}{(2s)!}\Bigg){\frac{n!}{(n-2s)!}\sum_{k\,=\,1}^{x-1}{\sum_{j\,=\,1}^{k}{j^{n-2s}}}},\]
and hence
\[x^n=x+\sum_{s\,=\,\frac{r}{2}, \ 1 \leq s \leq \frac{n}{2}}{\frac{2n!}{(2s)!(n-2s)!}\sum_{k\,=\,1}^{x-1}{\sum_{j\,=\,1}^{k}{j^{n-2s}}}}.\]
Finally, to determine how many iterations $s$ will take, that is, hows many times the nested summation $\frac{2n!}{(2s)!(n-2s)!}\sum_{k\,=\,1}^{x-1}{\sum_{j\,=\,1}^{k}{j^{n-2s}}}$ will be added, we note that $s$ essentially counts how many even numbers there are from $0$ to, and including, $n$. So, if $n$ is even, then there are $\frac{n}{2}$ even numbers $> 0 $ and $\leq n$; likewise, if $n$ is odd, then there are $\frac{n-1}{2}$ even numbers $ < n$. Therefore,
\[x^n=x+\sum_{s\,=\,1}^{m}{\frac{2n!}{(2s)!(n-2s)!}\sum_{k\,=\,1}^{x-1}{\sum_{j\,=\,1}^{k}{j^{n-2s}}}},\]
where $m= \begin{cases} 
      \frac{n}{2}& \text{if} \ n \textrm{ is even} \\
      \frac{n-1}{2}& \text{if} \ n \textrm{ is odd} 
\end{cases},$ 
and hence the theorem is proven. \hfill $\square$
\\
\\
\\
\indent The result of Theorem 2.1 provokes a search for even greater structures that underline the nature of integer powers and allow those to be represented by sums of smaller integral powers. In the next section, we consider what shall be called the ``bigger brother'' of (3).
\vspace{7cm}
\begin{center}{\textsc{3. Representation modulo 3}}\end{center}
\hspace{0.45cm} Looking back at Theorem 2.1, a natural thing to ponder is whether there is a similar identity in which $3$, and not $2$, plays a crucial role. Specifically, we think of something along the lines of
\begin{equation}x^n=x+\sum_{s\,=\,1}^{m}\frac{3n!}{(3s)!(n-3s)!}\Bigg(\sum_{k\,=\,1}^{x-1}\sum_{j\,=\,1}^{k}{\sum_{i\,=\,1}^{j}{i^{n-3s}}}\Bigg),\tag{$*$}\end{equation}
where $x,n,m \in \mathbb{Z^+}$, and three sums are summed $m$ times instead of two sums to match the pattern of $3$ used instead of $2$ everywhere. It is a rather tedious, but not an overly complicated calculation that one can perform to verify that, for the exponent $n=5$, the representative identity has the following appearance.
\[x^5=x+\sum_{s\,=\,1}^{\frac{5-2}{3}}\frac{3\cdot5!}{(3s)!(5-3s)!}\Bigg(\sum_{k\,=\,1}^{x-1}\sum_{j\,=\,1}^{k}{\sum_{i\,=\,1}^{j}{i^{5-3s}}}+\sum_{k\,=\,1}^{x-2}\sum_{j\,=\,1}^{k}{\sum_{i\,=\,1}^{j}{i^{5-3s}}}\Bigg).\]
On the other hand, for $n=4$, it is
\[x^4=x+\sum_{s\,=\,1}^{\frac{4-1}{3}}\frac{3\cdot4!}{(3s)!(4-3s)!}\Bigg(\sum_{k\,=\,1}^{x-1}\sum_{j\,=\,1}^{k}{\sum_{i\,=\,1}^{j}{i^{4-3s}}}+\sum_{k\,=\,1}^{x-2}\sum_{j\,=\,1}^{k}{\sum_{i\,=\,1}^{j}{i^{4-3s}}}\Bigg)+(x-1)x.\]
Notice that the above identities have two nested triple summations serving as the term of the summation that iterates over $s$ instead of just one, as it is in (3). Also, the fourth power case has a polynomial tailing at the end of the expression. The following conjecture generalizes representation modulo 3 and states every such tailing polynomial for every $n$.
\\ \\ \textbf{Conjecture 3.1} \textit{For all positive integers $x$ and $n$,}
\[x^n=x+\sum_{s\,=\,1}^{m}\frac{3n!}{(3s)!(n-3s)!}\Bigg(\sum_{k\,=\,1}^{x-1}\sum_{j\,=\,1}^{k}{\sum_{i\,=\,1}^{j}{i^{n-3s}}}+\sum_{k\,=\,1}^{x-2}\sum_{j\,=\,1}^{k}{\sum_{i\,=\,1}^{j}{(-1)^{s+1}i^{n-3s}}}\Bigg)+q_n(x)\]
\textit{with} $m= \begin{cases} 
      \frac{n}{3}& \text{if} \ n \equiv 0 \ \textrm{mod} \ 3 \\
      \frac{n-1}{3}& \text{if} \ n \equiv 1 \ \textrm{mod} \ 3 \\
      \frac{n-2}{3}& \text{if} \ n \equiv 2 \ \textrm{mod} \ 3
\end{cases}$, \ \textit{where} $q_n(x) = q_r(x)$ \textit{such that} {$n \equiv r \ \textrm{mod} \ 6$, $1 \leq r \leq 6$}. 
\\ \textit{The polynomials $q_r(x)$ for $r=1,2,3,4,5, 6$ are as follows.}
\[q_1=0,\] 
\[q_2=(x-1)x,\]
\[q_3=\frac{3}{2}(x-1)x,\]
\[q_4=(x-1)x,\]
\[q_5=0,\]
\[q_6=-\frac{1}{2}(x-1)x.\]
\begin{center}{\textsc{References}}\end{center}
\hypertarget{ref1}{$[1]$} \ \! \emph{Fermat's Little Theorem,} The Department of Mathematics and Computer Science, Emory University, \url{https://math.oxford.emory.edu/site/math125/fermatsLittleTheorem/}
\\ \\ \hypertarget{ref2}{$[2]$} \ \! \emph{two proofs of euclid's lemma,} Brooklyn College, City University of New York, \\ 
\url{http://www.sci.brooklyn.cuny.edu/~mate/misc/euclids_lemma.pdf}
\\ \\ \indent \textsc{\footnotesize{Department of Mathematics, California State University Long Beach, 1250 Bellflower Blvd, Long Beach, CA 90840, USA}}
\\ \indent \emph{Email address}: \href{mailto:nickgurin13@gmail.com}{\nolinkurl{nickgurin13@gmail.com}}
\end{document}